\documentclass[a4paper,11pt]{article}

\usepackage{MyLaTeX,MyArticle}
\input xy
\xyoption{all}

\title{On the Tensor Products of Modules for Dihedral $2$-Groups}
\author{David A.~Craven, University of Oxford}
\date{December 2007}

\begin{document}
\maketitle

%\begin{abstract}
\noindent Recall that an algebraic module is a $KG$-module that satisfies a polynomial with integer coefficients, with addition and multiplication given by direct sum and tensor product. In this article we prove that if $\Gamma$ is a component of the (stable) Auslander--Reiten quiver for a dihedral $2$-group consisting of non-periodic modules, then there is at most one algebraic module on $\Gamma$.
%\end{abstract}

\section{Introduction}

The only groups for which all indecomposable modules are `knowable' are those with cyclic, dihedral, semidihedral, and quaternion Sylow $p$-subgroups. The structure of the Green ring for groups with cyclic and $V_4$ Sylow $p$-subgroups are known, but no others have been determined. Of the remaining groups, the dihedral $2$-groups have the simplest module category but yet the tensor products of any two indecomposable modules are not known. In this article, we will prove that most non-periodic modules have a complicated tensor structure.

Following Alperin in \cite{alperin1976b}, we define a module to be \emph{algebraic} if it satisfies a polynomial with integer coefficients, where addition and multiplication are given by the direct sum and the tensor product. It is clear that a module $M$ is algebraic if and only if there are only finitely many isomorphism types of indecomposable summand in the collection of modules $M^{\otimes n}$ for all $n\geq 0$. Examples include all projective modules, more generally all trivial source modules, and all simple modules for $p$-soluble groups \cite{feit1980} and groups with abelian Sylow $2$-subgroups \cite{craven2007un2}.

\begin{thma}\label{mainthm} Let $\Gamma$ be a component of the (stable) Auslander--Reiten quiver $\Gamma_s(KG)$, where $G$ is a dihedral $2$-group. If $\Gamma$ contains non-periodic modules, then at most one module on $\Gamma$ is algebraic.
\end{thma}

If $G=V_4$, then a module is algebraic if and only if it is a sum of periodic modules and trivial modules. No such easy description is known for other dihedral $2$-groups, although there is the following conjecture.

\begin{conja} Let $M$ be a faithful, indecomposable module for a dihedral $2$-group $G$, and supposed that $M$ is not induced from a proper subgroup of $G$. Then $M$ is algebraic if and only if $M$ is periodic.
\end{conja}

This conjecture appears to be the first step in understanding the tensor structure of the modules for dihedral $2$-groups.

\section{Preliminaries on Module Theory}

We begin with the trivial results on algebraic modules.

\begin{lem}[{{\cite[Section II.5]{feit}}}] Let $M=M_1\oplus M_2$ be a $KG$-module, and suppose that $H_1\leq G\leq H_2$.
\begin{enumerate}
\item $M$ is algebraic if and only if $M_1$ and $M_2$ are algebraic.
\item The module $M_1\otimes M_2$ is algebraic.
\item The modules $M_1\res {H_1}$ and $M_1\ind {H_2}$ are algebraic.
\end{enumerate}
\end{lem}

An easy corollary of this lemma is that an indecomposable module is algebraic if and only if its source is.

\begin{thm}[{{\cite[Theorem 2.1]{bensoncarlson1986}}}]\label{bencarl} Let $G$ be a finite group and let $M$ and $N$ be absolutely indecomposable $KG$-modules, where $K$ is a field of characteristic $p$. Then $K|M\otimes N$ if and only if $p\nmid \dim M$ and $M\cong N^*$, in which case $K\oplus K$ is not a summand of $M\otimes N$. If $p\mid \dim M$, then every summand of $M\otimes N$ has dimension a multiple of $p$. 
\end{thm}

This result is obviously useful when dealing with tensor products. Finally, we include a lemma that will be required when dealing with the Auslander--Reiten quiver.

\begin{lem}[{{\cite[Proposition 4.12.10]{bensonvol1}}}]\label{splitrest} Let $M$ be an indecomposable module with vertex $Q$, and suppose that $H$ is a subgroup of $G$ not containing any conjugate of $Q$. Then the Auslander--Reiten sequence terminating in $M$ splits upon restriction to $H$.\end{lem}

\section{Indecomposable String Modules for Dihedral $2$-Groups}

In \cite{ringel1975}, Ringel classifies the indecomposable modules for the dihedral $2$-groups, and splits them into two collections: the \emph{string modules} and the \emph{band modules}. The band modules are all periodic, and so we will mostly ignore them in what follows. We assume that the reader is familiar with the construction of string modules, as given in \cite{ringel1975}, and we give one example to fix notation.

Write $\ms W$ for the set of strings of alternating $a^{\pm 1}$ and $b^{\pm 1}$. We call a symbol $a$ or $b$ a \emph{direct letter} and a symbol $a^{-1}$ or $b^{-1}$ an \emph{inverse letter}. Our modules are right modules, and so if $w=ab^{-1}aba^{-1}$, then the two matrices $\alpha$ and $\beta$ for $M(w)$ acting on the space $V$ with basis $\{v_1,\dots,v_6\}$ are given by
\[ \alpha=\begin{pmatrix}1&0&0&0&0&0\\1&1&0&0&0&0\\0&0&1&0&0&0\\0&0&1&1&0&0\\0&0&0&0&1&1\\0&0&0&0&0&1\end{pmatrix},\qquad \beta=\begin{pmatrix}1&0&0&0&0&0\\0&1&1&0&0&0\\0&0&1&0&0&0\\0&0&0&1&0&0\\0&0&0&1&1&0\\0&0&0&0&0&1\end{pmatrix}.\]
If $G=\Gen{x,y}{x^2=y^2=(xy)^{2q}=1}$, then let $M(w)$ denote the function $G\to \GL_n(2)$ defined by $x\mapsto \alpha$ and $y\mapsto \beta$. This will be a representation of the dihedral group $D_{4q}$ whenever no instance of $(ab)^q$, $(ba)^q$, $(a^{-1}b^{-1})^q$, or $(b^{-1}a^{-1})^q$ occurs. For the subset of $\ms W$ so defined, we use the symbol $\ms W_q$.

There are three important points to be made about the representations $M(w)$: firstly, they are always indecomposable representations; and secondly, $M(w)$ and $M(w')$ are isomorphic if and only if $w'=w$ or $w'=w^{-1}$. This latter point is crucial, and we will often blur the distinction between the words $w$ and $w^{-1}$. The last important point is that any odd-dimensional indecomposable module is a string module for some string of even length.

We need to briefly consider the band modules, to prove an easy fact about them, namely that for $M$ a band module, the modules $M\res{\gen x}$ and $M\res{\gen y}$ are both projective. We will not recall the definition of band modules here, but refer to \cite{ringel1975} for their construction. We will use the definition employed there.

\begin{lem}\label{resttogens} Let $M$ be an indecomposable $KG$-module.
\begin{enumerate}
\item If $M$ is odd-dimensional then $M\res{\gen x}$ and $M\res{\gen y}$ are both the sum of a trivial module and projective modules.
\item If $M$ is an even-dimensional string module then either $M\res{\gen x}$ is projective and $M\res{\gen y}$ is the direct sum of two copies of $K$ and a projective, or vice versa.
\item If $M$ is a band module, then both $M\res{\gen x}$ and $M\res{\gen y}$ are projective.
\end{enumerate}
\end{lem}
\begin{pf} Let $w$ be a word of even length $2n$, beginning with $a^{\pm 1}$ say, and let $v_i$ denote the standard basis, for $1\leq i\leq 2n+1$. Then the submodules of $M\res{\gen x}$ generated by $v_i$ and $v_{i+1}$ for $1\leq i<2n+1$ and $i$ odd form copies of projective modules, which therefore split off. Hence $M\res{\gen x}$ is the sum of $n$ projective modules and a trivial module. The same occurs for $M\res{\gen y}$, proving (i).

If $M$ is an even-dimensional string module then it is defined by a word $w$ of odd length $2n-1$, with first and last letters $a^{\pm 1}$ without loss of generality. Then $M\res{\gen y}$ has $n-1$ submodules $\gen{v_i,v_{i+1}}$ (for $i$ even) isomorphic with the projective indecomposable $K\gen y$-module, and two trivial submodules, $\gen{v_1}$ and $\gen{v_{2n}}$. Similarly, $\gen{v_i,v_{i+1}}$ is a projective submodule of $M\res{\gen x}$ for each odd $i$, and so $M\res{\gen x}$ is projective, proving (ii).

It remains to discuss the band modules. By cycling, we may assume that the word begins with $a$, and then we again see easily that the matrix corresponding to the action of $y$ on $M$ is a sum of projective modules, and this is true for any band module for a word beginning $a^{\pm 1}$. However, by cycling the word we find that $M$ is isomorphic with a band module for a word beginning $b^{\pm 1}$, and hence $M\res{\gen x}$ must also be projective, as required.
\end{pf}

Lemma \ref{resttogens}(i) allows us to define a group structure on the set of all odd-dimensional indecomposable modules, and in \cite{archer2008}, Archer studies this group, in particular proving Theorem \ref{mainthm} for this collection of modules. Therefore we need to understand even-dimensional string modules.

\begin{lem}\label{lemevenstring} Let $w,w'\in \ms W$ be words, and suppose that $\ell(w)=2n-1$ and $\ell(w')=2m-1$ are odd. Write $M=M(w)$ and $M'=M(w')$.
\begin{enumerate}
\item The word $w$ begins with $a^{\pm 1}$ if and only if it ends with $a^{\pm 1}$.
\item If $w$ begins with $a^{\pm 1}$, then the restriction $M\res{\gen x}$ is projective, and the restriction $M\res {\gen y}$ is the sum of a $2(m-1)$-dimensional projective module and a 2-dimensional trivial module.
\item If $w$ begins with $a^{\pm 1}$ and $w'$ begins with $b^{\pm 1}$, then $M\otimes M'$ contains no summands that are string modules.
\item If both $w$ and $w'$ begin with $a^{\pm 1}$, then $M\otimes M'$ contains exactly two even-dimensional string module summands.
\end{enumerate}
\end{lem}
\begin{pf} (i) is obvious, and (ii) easily follows from the construction of string modules, since the only place that a trivial summand can occur is at the beginning or end of a word. The proof of (iii) comes from the fact that if $M\otimes M'$ contains a string module, there must be a trivial summand of either $(M\otimes M')\res{\gen x}$ or $(M\otimes M')\res{\gen y}$, which is impossible since both $M\res{\gen x}$ and $M'\res{\gen y}$ are projective. The proof of (iv) is similar: if $M$ and $M'$ both begin with $a^{\pm 1}$, then both $M\res{\gen y}$ and $M'\res{\gen y}$ contain two trivial summands, proving that $(M\otimes M')\res{\gen y}$ contains four trivial summands. Since band modules restrict to projective modules, and no odd-dimensional summand can occur by Theorem \ref{bencarl}, the tensor product must contain two even-dimensional string modules as summands.\end{pf}

Write $z$ for the non-trivial central element, and write $X=\gen{x,z}$ and $Y=\gen{y,z}$. By the Alperin--Evens theorem \cite{alperinevens1981}, if $M$ is a non-periodic module, either $M\res X$ or $M\res Y$ is non-periodic (since $X$ and $Y$ are representatives for the two conjugacy classes of $V_4$ subgroup).

Suppose, without loss of generality, that $w$ begins with $a^\ep$, so that $M\res {\gen x}$ is projective and $M\res {\gen y}$ is non-projective. Since $\gen x$ has index 2 in $X$, it must be true that $M\res X$ is periodic, and so $M\res Y$ is non-periodic. It is well-known (and a consequence of the construction of the string modules) that the only non-periodic modules for $V_4$ are the Heller translates of the trivial module. It can easily be seen that $M\res Y$ must be the sum of two odd-dimensional modules $\Omega^r(K)\oplus \Omega^s(K)$ and periodic modules. We call the \emph{signature} of the module $M$ the object $[r,s]$, where
\[ \Omega^r(K)\oplus \Omega^s(K)|M_{(i,j)}\res Y.\]

We now need to understand the Auslander--Reiten quiver. In order to describe the action of $\Omega^2$ on string modules effectively, we introduce two operations, $L_q$ and $R_q$, on the set of all words $\ms W_q$. Write $A=(ab)^{q-1}a$ and $B=(ba)^{q-1}b$. The operator $L_q$ is defined by adding or removing a string at the start of the word $w$, and $R_q$ is the same but at the end of the word.

If the word $w$ starts with $Ab^{-1}$ or $Ba^{-1}$, then $wL_q$ is $w$ with this portion removed. If neither of these are present, then we add either $A^{-1}b$ or $B^{-1}a$ to $w$ to get $wL_q$, whichever gives an element of $\ms W_q$. Similarly, if $w$ ends with $aB^{-1}$ or $bA^{-1}$, then $wR_q$ is $w$ with this portion removed. If neither of these are present, then we add either $a^{-1}B$ or $b^{-1}A$ to $w$ to get $wR_q$, whichever gives a word in $\ms W_q$. The operators $L_q$ and $R_q$ commute, and are bijections on $\ms W_q$.

The double Heller operator $\Omega^2$ is given by
\[ \Omega^2(M(w))=M(wL_qR_q),\]
and the almost-split sequences on string modules are given by
\[ 0\to M(wL_qR_q)\to M(wL_q)\oplus M(wR_q)\to M(w)\to 0,\]
unless $w=AB^{-1}$, in which case the almost-split sequence is
\[ 0\to M(wL_qR_q)\to M(wL_q)\oplus M(wR_q)\oplus KG\to M(w)\to 0,\]
where $KG$ denotes the projective indecomposable module $KG$, viewed as a module over itself. This describes the Auslander--Reiten quiver, and it looks as follows.
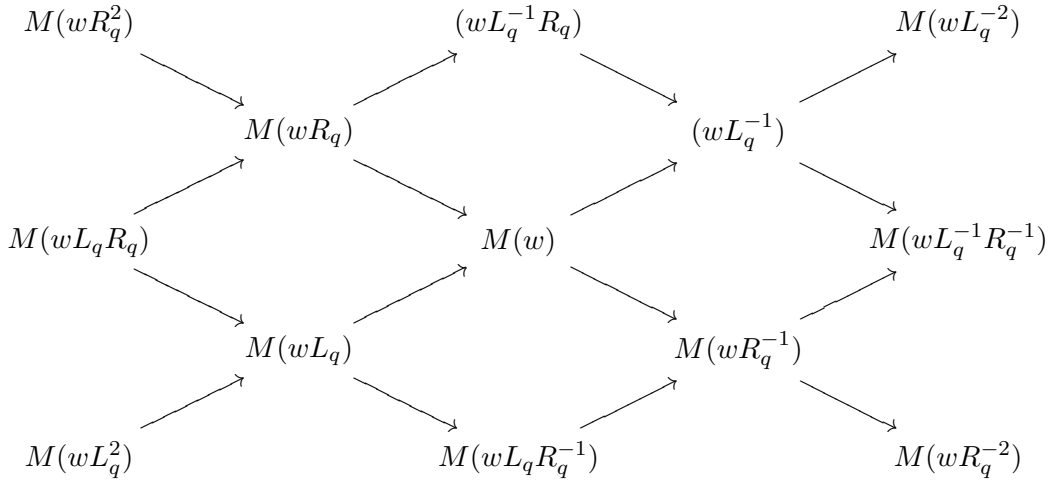
\[ \begin{diagram}\node{M(wR_q^2)}\arrow{se}\node{} \node{(wL_q^{-1}R_q)}\arrow{se}\node{} \node{M(wL_q^{-2})}
\\ \node{}\node{M(wR_q)}\arrow{se}\arrow{ne}\node{} \node{(wL_q^{-1})}\arrow{se}\arrow{ne}\node{}
\\ \node{M(wL_qR_q)}\arrow{se}\arrow{ne} \node{}\node{M(w)}\arrow{se}\arrow{ne} \node{}\node{M(wL_q^{-1}R_q^{-1})}
\\ \node{}\node{M(wL_q)}\arrow{se}\arrow{ne} \node{}\node{M(wR_q^{-1})}\arrow{se}\arrow{ne}\node{}
\\ \node{M(wL_q^2)}\arrow{ne} \node{}\node{M(wL_qR_q^{-1})}\arrow{ne} \node{}\node{M(wR_q^{-2})}
\end{diagram}\]
[In this diagram, the $\Omega^2$ operation is a functor moving from right to left, the map $M(w)\mapsto M(wL_q)$ is a function moving down and to the left, and the map $M(w)\mapsto M(wR_q)$ moves up and to the left.]

Considering a component $\Gamma$ of the $\Gamma_s(KG)$, we will abuse notation slightly, and also refer to the signature of a vertex, as well as the signature of a module.

\section{Algebraicity of Modules}

As we have mentioned, in \cite[Theorem 3.4]{archer2008}, Archer proves that there are no non-trivial, indecomposable algebraic modules of odd dimension. Thus Theorem \ref{mainthm} reduces to proving the result for components of $\Gamma_s(KG)$ containing even-dimensional string modules. It suffices to prove the following result.

\begin{thm}\label{mainthm2} Let $\Gamma$ be a component of $\Gamma_s(KG)$ containing non-periodic modules of even dimension. Then there is at most a single module on $\Gamma$ with signature $[0,0]$.
\end{thm}

We will prove Theorem \ref{mainthm2} in a sequence of lemmas. We begin with the following observation.

\begin{lem}\label{fixpoints} Let $G=V_4$, and let $x$ be a non-identity element of $G$. Let $i$ be a non-positive integer, and let $M=\Omega^i(K)$. Then the $G$-fixed points of $M$ are equal to the $x$-fixed points of $M$.
\end{lem}
\begin{pf} It is easy to see that the socle of $M$ is of dimension $i+1$. We simply note that $M\res{\gen x}$ is the sum of $K$ and $i$ copies of the free module, and so its socle has dimension $i+1$ also. Thus the lemma must hold.\end{pf}

Using this lemma, we can prove a crucial result about the summands of $M(w)\res Y$ under a certain condition on $w$.

\begin{lem}\label{fixpointexists} Suppose that $M=M(w)$ is an even-dimensional string module, and suppose that $w$ begins with $a^{-1}$ or ends with $a$. Finally, suppose that the odd-dimensional summands of $M\res Y$ are isomorphic with $\Omega^i(K)$ and $\Omega^j(K)$, where both $i$ and $j$ are non-positive. Then (at least) one of $i$ and $j$ is 0.
\end{lem}
\begin{pf} Note that, since $w$ begins with an inverse letter, the subspace $U=\gen{v_i:i\geq 2}$ is a $G$-submodule of $M$ (where the $v_i$ are the standard basis used in the construction of the string modules). Thus if there exists a $Y$-fixed point
\[ V= v_1+\sum_{i\in I} v_i,\]
then $\gen V$ is a summand of $M\res Y$ isomorphic with $K$, as required. Let $N_1$ and $N_2$ denote the two odd-dimensional summands of $M\res Y$. By Lemma \ref{fixpoints}, it suffices to show that there is such a point $V$ fixed by $y$ lying inside one of the $N_i$.

We will now calculate the possibilities for a trivial summand of $M\res{\gen y}$. Since $\gen{v_2,\dots,v_{n-1}}\res{\gen y}$ (where $\dim M=n$) is a free module, if $\alpha=\sum_{j\in J} v_j$ is a fixed point of $M\res{\gen y}$ with a complement, then either $1$ or $n$ lies in $J$. Since $M\res{\gen y}$ contains a 2-dimensional trivial module, we easily see that the fixed points with complements are given by
\[ v_1+\sum_{j\in J} v_j, \quad v_n+\sum_{j\in J} v_j,\quad v_1+v_n+\sum_{j\in J} v_j,\]
where $J\subs \{2,\dots,n-1\}$. Hence for some suitable choice of $I$, the point $V$ given above is a $y$-fixed point, as required.

If $w$ ends with $a$, then $w^{-1}$ begins with $a^{-1}$. Since $M(w)=M(w^{-1})$, we get the result.
\end{pf}

As a remark, by taking duals, one sees that if $M=M(w)$ and $w$ begins with $a$ or ends with $a^{-1}$, and the odd-dimensional summands of $M\res Y$ are isomorphic with $\Omega^i(K)$ and $\Omega^j(K)$ for $i,j\geq 0$, then (at least) one of $i$ and $j$ is 0.

To provide the proof of Theorem \ref{mainthm2}, we must analyze the components of the Auslander--Reiten quiver consisting of non-periodic, even-dimensional string modules. To do this, let $M$ denote such an indecomposable module, and suppose without loss of generality that $M=M(w)$ where $w$ begins with $a^{\pm 1}$. Denote by $\Gamma$ the component of $\Gamma_s(KG)$ on which $M$ lies.

We will co-ordinatize the component $\Gamma$: write $(0,0)$ for the co-ordinates of the vertex corresponding to $M(w)$, and $(i,j)$ for the vertex corresponding to $M(wL_q^iR_q^j)$. Then the portion of $\Gamma$ around the module $M$ is co-ordinatized as follows.
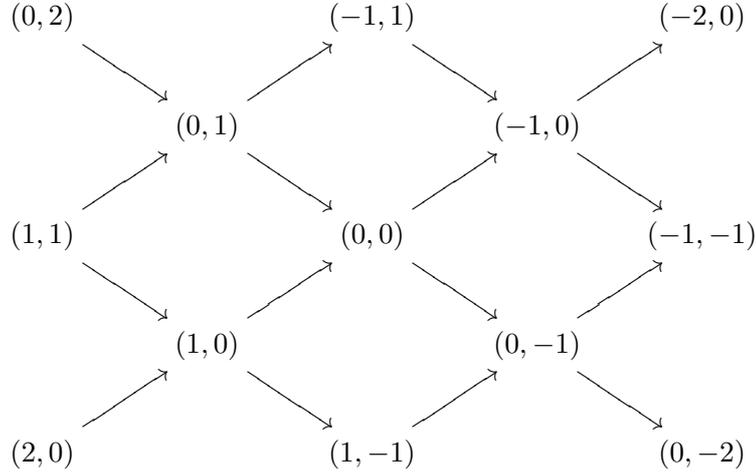
\[ \begin{diagram}\node{(0,2)}\arrow{se}\node{} \node{(-1,1)}\arrow{se}\node{} \node{(-2,0)}
\\ \node{}\node{(0,1)}\arrow{se}\arrow{ne}\node{} \node{(-1,0)}\arrow{se}\arrow{ne}\node{}
\\ \node{(1,1)}\arrow{se}\arrow{ne} \node{}\node{(0,0)}\arrow{se}\arrow{ne} \node{}\node{(-1,-1)}
\\ \node{}\node{(1,0)}\arrow{se}\arrow{ne} \node{}\node{(0,-1)}\arrow{se}\arrow{ne}\node{}
\\ \node{(2,0)}\arrow{ne} \node{}\node{(1,-1)}\arrow{ne} \node{}\node{(0,-2)}
\end{diagram}\]

We get a `diamond rule' for the diamonds of the Auslander--Reiten quiver using Lemma \ref{splitrest}, so that if $M_{(i,j)}$ does not have vertex contained within $Y$, then
\[ M_{(i,j)}\res Y\oplus M_{(i+1,j+1)}\res Y=M_{(i,j+1)}\res Y\oplus M_{(i+1,j)}\res Y.\]

Suppose that no module on $\Gamma$ has vertex $Y$. (Since every proper subgroup of $Y$ is cyclic, if $N$ is a non-periodic indecomposable module with vertex contained within $Y$, it has vertex $Y$.) If the signatures are known for two adjacent rows of $\Gamma$, then they can be calculated for all rows, using the diamond rule. Since two rows (say rows $\alpha$ and $\alpha+1$) are completely known, the rows $\alpha+2$ and $\alpha-1$ can be calculated, since every point on either of those rows lies on a diamond whose other three corners lie in the rows $\alpha$ and $\alpha+1$. This process can be iterated to get the signatures for all rows.

This information makes the proof of the next proposition possible.

\begin{prop}\label{trichotARquiver} Let $M=M_{(0,0)}$ be a non-periodic, even-dimensional string module, and suppose that $M$ is algebraic. Suppose in addition that the component $\Gamma$ of $\Gamma_s(KG)$ containing $M$ contains no module with vertex $Y$. Let $M_{(i,j)}$ denote the indecomposable module $M(wL_q^i R_q^j)$. Write $[r,s]$ for the signature of $(i,j)$. Then exactly one of the following three possibilities occurs:
\begin{enumerate}
\item The signature of $(i,j)$ is $[2i,2j]$ (or $[2j,2i]$);
\item The signature of $(i,j)$ is $[2i,2i]$; and
\item The signature of $(i,j)$ is $[2j,2j]$.
\end{enumerate}
\end{prop}
\begin{pf}Firstly, we note that all three potential signatures satisfy the diamond rule that the sum of the signatures of $(i,j)$ and $(i-1,j-1)$ is equal to the sum of the signatures of $(i-1,j)$ and $(i,j-1)$. We need to check that these three possibilities are the only ones, and by the remarks before the proposition it suffices to check that these are the only three possibilities for the two rows with vertices $(i,i)$ and $(i,i+1)$ in the Auslander--Reiten quiver.

Since the signature of $(0,0)$ is $[0,0]$, the signature of $(i,i)$ must be $[2i,2i]$, since
\[ M_{(i,i)}=\Omega^{2i}(M_{(0,0)}).\]

Since no module on $\Gamma$ has vertex contained within $Y$, the diamond rule for the diamond containing $(0,0)$ and $(1,1)$ becomes
\[ M_{(0,0)}\res Y\oplus M_{(1,1)}\res Y=M_{(0,1)}\res Y\oplus M_{(1,0)}\res Y.\]
The signatures of $(0,0)$ and $(1,1)$ are $[0,0]$ and $[2,2]$ respectively, and so the signature of $(0,1)$ is one of $[0,2]$ (or equivalently $[2,0]$), $[0,0]$ or $[2,2]$. Thus the signatures of $(i,i+1)$ are one of $[2i,2i+2]$, $[2i,2i]$ or $[2i+2,2i+2]$, which correspond to (i), (ii) and (iii) respectively in the proposition.
\end{pf}

In fact, the same result holds for the two components containing non-periodic modules with vertex $Y$, but it requires more work.

Let $M$ be an indecomposable module with vertex $Y$. If $M$ is non-periodic, then the source $S$ of $M$ must also be non-periodic. Thus $S=\Omega^i(K)$ for some $i\in \Z$. Therefore the modules $\Omega^i(K_Y)\ind G$ (where $K_Y$ denotes the trivial module for $Y$) are the only non-periodic indecomposable modules with vertex $Y$. The module $(K_Y)\ind G$ is algebraic, whereas all others are not.

We begin by considering the component containing $M_{(0,0)}=\Omega(K_Y)\ind G$. This cannot contain algebraic modules, because it can have no vertex with signature $[0,0]$. To see this, notice firstly that the signature of $(0,0)$ is $[1,1]$. We analyze the diamond with bottom vertex $(0,0)$: write $[r,s]$ for the signature of the top vertex, namely $(-1,1)$, and write $[p,q]$ for the signature of the vertex $(0,1)$ on the left of the diamond. Then the diamond rule gives
\[ [1,1]\union[r,s]=[p,q]\union[p-2,q-2],\]
and we see that $p$, $q$, $r$ and $s$ are all odd. Thus all signatures of vertices $(i,i+1)$ (i.e., the row above that containing $M_{(0,0)}$) are a pair of odd numbers. Since all diamonds not involving those modules with vertex $Y$ obey the diamond rule, we see that all modules above the horizontal line containing $M_{(0,0)}$ have signature a pair of odd numbers. The same analysis holds for the lower half of the quiver, and so our claim holds.

The other component with modules of vertex $Y$, namely that containing $M_{(0,0)}=(K_Y)\ind G$, does contain an algebraic module. Suppose that the signatures of the vertices on the horizontal line containing $(0,0)$, and those on the lines directly above and below this are known. (Thus the signatures for all vertices $(i,i)$, $(i+1,i)$ and $(i-1,i)$ are known.) Then we claim that the signatures for all vertices can be deduced. This is true for the same reason as before, since all diamonds containing at most one point from the line of vertices $(i,i)$ obey the diamond rule.

This will enable us to prove the next proposition easily.

\begin{prop}\label{trichotARquiver2} Let $M=M_{(0,0)}$ be the module $K_Y\ind G$, where $K_Y$ denotes the trivial module for $Y$. Let $M_{(i,j)}$ denote the indecomposable module $M(wL_q^i R_q^j)$. Write $[r,s]$ for the signature of $(i,j)$. Then exactly one of the following three possibilities occurs:
\begin{enumerate}
\item The signature of $(i,j)$ is $[2i,2j]$ (or $[2j,2i]$);
\item The signature of $(i,j)$ is $[2i,2i]$; and
\item The signature of $(i,j)$ is $[2j,2j]$.
\end{enumerate}
\end{prop}
\begin{pf} Firstly note that the three signature patterns obey the diamond rule everywhere, so they certainly obey it for those diamonds that split upon restriction to $Y$. Thus we need only show that these three possibilities are the only ones. By the preceding remarks, it suffices to show this for the horizontal lines containing the vertices $(i,i)$, $(i,i-1)$ and $(i-1,i)$.

We analyze the diamond with bottom vertex $(0,0)$: write $[r,s]$ for the signature of the top vertex, namely $(-1,1)$, and write $[p,q]$ for the signature of the vertex $(0,1)$ on the left of the diamond. Then the diamond rule gives
\[ [0,0]\union[r,s]=[p,q]\union[p-2,q-2],\]
and so $p$ and $q$ are either both $0$, both $2$, or one is $0$ and one is $2$. In any case, this uniquely determines all modules on the horizontal line containing the vertex $(0,1)$, and they are as claimed in the proposition. We need to determine the signatures of the vertices $(i,i-1)$ from these.

Suppose that the signature of $M_{(0,1)}$ is $[0,0]$. Then the dual of $M_{(0,1)}$ must also have signature $[0,0]$. The almost-split sequence terminating in $M_{(0,0)}$ is given by
\[ 0\to M_{(1,1)}\to M_{(0,1)}\oplus M_{(1,0)}\to M_{(0,0)}\to 0,\]
and since $M_{(0,0)}$ is self-dual, the dual of this sequence is the (almost-split) sequence
\[ 0\to M_{(0,0)}\to M_{(0,-1)}\oplus M_{(-1,0)}\to M_{(-1,-1)}\to 0.\]
Thus either $M_{(0,1)}^*=M_{(0,-1)}$ or $M_{(0,1)}^*=M_{(-1,0)}$. However, the second possibility cannot occur, since we know that the signature of $(-1,0)$ is $[-2,-2]$, and thus
\[ M_{(0,1)}^*=M_{(0,-1)}.\]
Hence the signature of $(0,-1)$ is $[0,0]$, and we have proved that the three lines containing the vertices $(i,i)$, $(i,i-1)$ and $(i-1,i)$ have signatures obeying possibility (ii).

Now suppose that the signature of $M_{(0,1)}$ is $[2,2]$. Then $M_{(0,1)}^*\not\cong M_{(-1,0)}$ since the signature of $M_{(-1,0)}$ is $[0,0]$. Thus we again have 
\[ M_{(0,1)}^*=M_{(0,-1)}.\]
Since the signature of $(0,1)$ is $[2,2]$, the signature of $(0,-1)$ is $[-2,-2]$, and so we have proved that the three lines containing the vertices $(i,i)$, $(i,i-1)$ and $(i-1,i)$ have signatures obeying possibility (iii).

Finally, suppose that the signature of $(0,1)$ is $[0,2]$. If the signature of $M_{(0,-1)}$ is not $[0,-2]$, then its dual would have to be $M_{(0,1)}$, by the same reasoning as the previous two paragraphs. However, this is not possible, and so we have proved that the three lines containing the vertices $(i,i)$, $(i,i-1)$ and $(i-1,i)$ have signatures obeying possibility (i).
\end{pf}

In the first case of Propositions \ref{trichotARquiver} and \ref{trichotARquiver2}, there is a unique vertex on $\Gamma$ with signature $[0,0]$, namely the vertex $(0,0)$, and so $M$ is indeed the unique algebraic module on $\Gamma$. This is in accordance with Theorem \ref{mainthm2}.

In the second case, $K\oplus K|M(wL_q^i)\res Y$ for all $i\in\Z$, and 
\[\Omega^{-2}(K)\oplus \Omega^{-2}(K)|M(wL_q^iR_q^{-1})\res Y.\]
If $i$ is a suitably large negative number, then $wL_q^iR_q^{-1}$ begins with $a^{-1}$. This yields a contradiction, since by Lemma \ref{fixpointexists}, $K$ must be a summand of $M(wL_q^iR_q^{-1})\res Y$.

In the third case, $K\oplus K|M(wR_q^i)\res Y$ for all $i\in\Z$, and so
\[ \Omega^2(K)\oplus \Omega^2(K)|M(wL_q^{-1}R_q^i)\res Y.\]
If $i$ is a suitably large negative number, then $wL_q^{-1}R_q^i$ ends with $a^{-1}$. This yields a contradiction, since by Lemma \ref{fixpointexists}, $K$ must be a summand of $M(wL_q^{-1}R_q^i)\res Y$.

Thus in Propositions \ref{trichotARquiver} and \ref{trichotARquiver2} only the first possibility can occur, and so Theorem \ref{mainthm2} is proved.

\end{document}